\documentstyle[12pt]{article}
\newcommand{\C}{\mbox{l\hspace{-.47em}C}}

\begin{document}

\title{INVARIANCE QUANTUM GROUPS OF THE DEFORMED OSCILLATOR ALGEBRA}
\author{Jacqueline BERTRAND and Mich\`{e}le IRAC-ASTAUD \\ 
Laboratoire de Physique Th\'{e}orique et Math\'{e}matique\\
Universit\'{e} Paris VII\\2 place Jussieu F-75251 Paris Cedex 05, FRANCE}
\date{}
\maketitle
submitted to Journal of Physics A : Mathematical and General.

\begin{abstract}
A differential calculus is set up on a deformation of the oscillator algebra. 
It is  uniquely determined by the requirement of invariance 
under a seven-dimensional quantum group. The quantum space and its 
associated differential calculus are also shown to be invariant under a 
nine generator quantum group containing the previous one.
\end{abstract}

\section{Introduction}
Several kinds of differential calculi have been introduced on quantum 
groups and quantum spaces and widely studied  (\cite{Woro}-\cite{wz}). 
However there is one important case that deserves a special 
treatment because of its importance in physics, namely 
the case of the oscillator (or Weyl-Heisenberg) algebra. Indeed, in its usual 
form it represents the algebra of observables in quantum mechanics. 
After deformation, it is still an algebra of observables but for a 
different quantization.   \par
The problem of constructing a differential calculus on such an 
algebra can be tackled in different ways. 
Here we privilege the invariance approach. 
Given the quantum space of variables (or observables), we 
first construct the set of seven element quantum matrices    preserving 
that space (section 2). The space of differentials is then determined 
in two independent ways: either by postulating the existence of 
a $R$-matrix (section 3) or by constructing directly a new invariant 
space (section 4). The uniqueness of the result emphasizes the  
power of the approach based on covariance. Finally, in section 5, the
quantum group preserving simultaneously the 
spaces of variables and differentials
is explicitly described. Moreover, a larger quantum group with nine  
generators is also shown to preserve the same spaces. Both 
groups can be endowed with a structure of Hopf algebra.\par

\section{The quantum space and its   invariance algebra}
The Weyl-Heisenberg algebra written in homogeneous form is the free 
associative algebra generated by three operators $x^i$ satisfying the 
following quadratic relations :
\begin{equation}
(R) \, \left\{
\begin{array}{l}
x^1 x^2 - x^2 x^1 - s(x^3)^2 =0\\
x^1 x^3 = x^3 x^1   \\
x^2 x^3 = x^3 x^2 
\end{array}
\right.
\label{1}
\end{equation}
This algebra will be denoted by $\C <x>/R$. It is invariant under the 
seven-parameter Lie subgroup 
$G$ of $GL(3)$ consisting of matrices $T$ such that 
$T^3_{\, 1}   =T^3_{\, 2}  =0$ and 
$T^1_{\, 1}  T^2_{\, 2} -T^1_{\, 2} T^2_{\, 1} =(T^3_{\, 3} )^2$.
\par
Deforming this algebra, 
we consider 
the quantum space obtained when 
  relations (\ref{1}) are replaced by
\begin{equation}
(R_{xx}) \hspace{5mm}
\left\{
\begin{array}{rcl}
x^1 x^2 - q \, x^2 x^1 -s (x^3)^2 & = & 0 \\
x^1x^3- u \, x^3 x^1 &  = & 0 \\
x^2x^3-u^{-1}\, x^3 x^2 & = & 0
\end{array}
\right.
\label{3}
\end{equation}
Let $V$ be the vector space of column vectors $X$ 
with elements $x^i$. The matrix $T$ now has non commuting 
elements and acts on $X$ according to the definition
 \cite{rtf} : 
\begin{equation}
\delta : X \longrightarrow \delta(X)= T\otimes X
\label{7}
\end{equation}
Relations
$R_{xx}$ defined by (\ref{3}) 
will be preserved by this action provided the 
following commutation relations hold:
\begin{equation}
\begin{array}{lll}
 T^1_{\, 1} T^3_{\, 3} =T^3_{\, 3} T^1_{\, 1}  & \hspace{1cm} &  T^2_{\, 1} T^3_{\, 3}  = u^{-2}\, T^3_{\, 3} T^2_{\, 1}  \\
 T^1_{\, 2} T^3_{\, 3} =u^2\, T^3_{\, 3} T^1_{\, 2}  & \hspace{1cm} & T^2_{\, 2} T^3_{\, 3}  = T^3_{\, 3} T^2_{\, 2}  \\
T^1_{\, 3} T^3_{\, 3}  = u \, T^3_{\, 3} T^1_{\, 3}   & \hspace{1cm} & T^2_{\, 3} T^3_{\, 3}  = u^{-1} \, T^3_{\, 3} T^2_{\, 3} \\
 T^1_{\, 1} T^2_{\, 1}   =  q \, T^2_{\, 1} T^1_{\, 1}  & \hspace{1cm} &   T^1_{\, 2} T^2_{\, 2}   =  q \, T^2_{\, 2} T^1_{\, 2} 
\end{array}
\label{8}
\end{equation}
\par
\vspace*{-6mm}
\begin{eqnarray}
u \, T^1_{\, 1} T^2_{\, 3}  - q \, T^2_{\, 3} T^1_{\, 1}  & = &  q u \, T^2_{\, 1} T^1_{\, 3}     - T^1_{\, 3} T^2_{\, 1}   
\label{81} \\  
  T^1_{\, 2} T^2_{\, 3}  - qu \, T^2_{\, 3} T^1_{\, 2}  & =  & q \, T^2_{\, 2} T^1_{\, 3}  - u \, T^1_{\, 3} T^2_{\, 2}        
\label{82} \\  
  T^1_{\, 1} T^2_{\, 2}  -T^2_{\, 2} T^1_{\, 1}  & = & q \, T^2_{\, 1} T^1_{\, 2}  - q^{-1} \, T^1_{\, 2} T^2_{\, 1}          
\label{83} \\ 
  (T^1_{\, 1} T^2_{\, 2}  - q\, T^2_{\, 1} T^1_{\, 2} )s & = & s\, (T^3_{\, 3} )^2 -T^1_{\, 3} T^2_{\, 3}  + q \, T^2_{\, 3} T^1_{\, 3}     
\label{84} 
\end{eqnarray}
In that case, the action $\delta $ becomes defined as a 
mapping of $\C <x>/R_{xx}$ onto itself.  

\section{ $R$-matrix and  invariant differential calculus }

In the present case, 
a differential calculus on $\C <x>/R_{xx}$ that is invariant under the 
action of  $T$  can be constructed if there exists a matrix 
$\widehat{R}$ with the following properties \cite{rtf} :
\par
 - \noindent $\widehat{R}$ is defined by the relations
\begin{equation}
\widehat{R}^{ji}_{kl} \, 
T^k_{\, m} T^l_{\, n} = T^j_{\, l} T^i_{\, k} \, \widehat{R}^{lk}_{mn}
\label{66}
\end{equation}

 - $\widehat{R}$ has two eigenspaces,  $V_1$ and $V_2$
that can  be identified with the variables and the 
 one-forms quantum spaces respectively. Space $V_1$   
has dimension six and is determined by the relations ($R_{xx}$) 
given in (\ref{3}). 
\par
The determination of $\widehat{R}$  is performed by assuming that 
relations (\ref{8}-\ref{84})  can be cast into the form (\ref{66}) and 
solving the corresponding equations. 
In addition, we impose two natural requirements:
\begin{description}
\item[{\it (i)}]
The determinant of $T$ is different from zero so that the 
set of matrices $T$ can be made into a Hopf algebra.
\item[{\it (ii)}]
The ordering of monomials such as $T^i_{\, j} T^k_{\, l} T^m_{\, n}$ 
is independent of the procedure used, when the associativity of the 
algebra $C<T>$ is taken into account.
\end{description}
With these conditions, it is found that a $\widehat{R}$-matrix  
exists only if $q=u^2$ and it is given by:
\begin{equation}
\widehat{R}
=\left(\begin{array}{c c c c c c c c c }
1&0&0&0&0&0&0&0&0\\
0&0&0&u^2&0&0&0&0&s\\
0&0&0&0&0&0&u&0&0\\
0&u^{-2}&0&0&0&0&0&0&-s/u^2\\
0&0&0&0&1&0&0&0&0\\
0&0&0&0&0&0&0&1/u&0\\
0&0&1/u&0&0&0&0&0&0\\
0&0&0&0&0&u&0&0&0\\
0&0&0&0&0&0&0&0&1\\
\end{array}\right)  
\label{25}
\end{equation}
It can be shown that 
this matrix 
$\widehat{R}$
 is equal to its inverse
and    verifies the Yang-Baxter equation:
\begin{equation}
(\widehat{R} \otimes 1)(1\otimes \widehat{R})(\widehat{R}\otimes 1)=(1\otimes \widehat{R})(\widehat{R}\otimes 1)(1\otimes \widehat{R})
\label{yb1}
\end{equation}
where $1$ is the unit matrix of $GL(3)$.
The matrix $\widehat{R}$ has two eigenspaces that correspond to the 
variables quantum space defined by the relations $R_{xx}$ and to the 
one-forms quantum space defined by :
\begin{equation} 
R_{\xi\xi} : \,
\left
\{\begin{array}{cclccl}
(\xi^1)^2 & = & 0,\quad & (\xi^2)^2 & = & 0\quad \\
(\xi^3)^2 & = & 0, \quad &\xi^2 \xi^1& = &-  u^{-2} \, \xi^1\xi^2 \\  
\xi^1 \xi^3&=&- u \, \xi^3 \xi^1, \quad 
& \xi^2 \xi^3 &= & - u^{-1} \, \xi^3 \xi^2 
\end{array}
\right.
\label{xi1}
\end{equation}
{\bf Result}:
If    the existence of 
a 
$\widehat{R}$-matrix is assumed and if conditions {\it (i)} and {\it (ii)} are 
satisfied, then  
an invariant differential calculus can be set up on $C<x>/R_{xx}$ if and only 
if $q = u^2$ . 
\par
In the next section, the same result is obtained without 
assuming the existence of a $\widehat{R}$-matrix.

\section{Invariant exterior algebra}
Consider a vector $\Xi  \in V$ with components $\xi^i$, $i=1,2,3$  
and form $V\otimes V$. 
The space of invariant forms can be constructed directly  
as an invariant 
subspace $V_2$ of $V\otimes V$ that is supposed to be of dimension three or less.
It will be done by 
finding the possible invariant relations  
between the $\xi^i \xi^j$. 
\par
An important tool for this derivation  
is the introduction of a degree $d^\circ$ on $C<T>$,  
i.e. of a homomorphism from $C<T>$ into $Z$. Explicitly, the 
degree operation $d^\circ$ 
associates with each element $T^i_{\, j}$  the power of $u$ present 
in the commutation relations (\ref{8}) of that  element  with $T^3_{\, 3}$. 
Assuming that  $u \ne 1$, we find the degrees of all the elements 
of $T$:
\begin{equation}
\begin{array}{llll}
d^\circ(T^1_{\, 1} )&=d^\circ(T^2_{\, 2} )&=d^\circ(T^3_{\, 3} )&=0\\
d^\circ(T^1_{\, 2} )&=2&&\\
d^\circ(T^2_{\, 1} )&=-2&&\\
d^\circ(T^1_{\, 3} )&=1&&\\
d^\circ(T^2_{\, 3} )&=-1&&
\end{array}
\end{equation}
The degrees of all the monomials in $C<T>$ can then be deduced.
\par
Consider a quadratic relation $R_1=0$ between the $\xi$.   
After action of the homomorphism $\delta \otimes \delta$ defined
in (\ref{7}), a new relation  $(T\otimes T)\otimes R_1=0$ is obtained 
which contains monomials  $T^i_{\, j} T^k_{\, l}$  of 
known degrees. Invariance of relation $R_1 = 0$ then implies  
that terms of different degrees vanish separately. 
Applying this remark systematically allows to
find all possible invariant relations as will now be sketched.
\par
Let  $\xi^{\prime i}\xi^{\prime j}$  denote the components   
of the vector 
$\Xi^\prime\otimes \Xi^\prime  \equiv (T\otimes T)\otimes ( \Xi\otimes \Xi)$. 
The monomial $(\xi^{\prime 1})^2$ 
is the only one  containing a term of degree 4. 
Hence it cannot be involved in any quadratic relation 
except $(\xi^1)^2=0$. The same reasoning with  $(\xi^{\prime 2} )^2$
and degree (-4) leads to $(\xi^2)^2=0$. In the remaining transformed 
relations,  terms of degree $(+2)$ (resp $(-2)$) come from monomials 
$\xi'^1\xi'^3$ and $\xi'^3\xi'^1$  (resp $\xi'^2\xi'^3$ and $\xi'^3\xi'^2$). 
This implies that  $\xi^1\xi^3$ and $\xi^2\xi^3$  are independent monomials.
To generate $V_2$, we need only   a third independent  monomial 
which can be chosen  as  
$\xi^1\xi^2, \xi^2\xi^1$  or $(\xi^3)^2$.
Choosing $\xi^1\xi^2$, we can write the  general form of the 
invariant relations as:
\begin{equation}
\begin{array}{lll}
& (\xi^1)^2&=0\\
& (\xi^2)^2&=0\\
& (\xi^3)^2&= k  \; \xi^1 \xi^2\\
& \xi^2 \xi^1 &= \xi \; \xi^1 \xi^2 \\
& \xi^3 \xi^1 &= \lambda \; \xi^1 \xi^3 + \lambda_{12} \; \xi^1 \xi^2\\
& \xi^3 \xi^2 &= \mu \; \xi^2 \xi^3+\mu_{12} \; \xi^1 \xi^2\\
\end{array}
\label{13}
\end{equation}
The conditions of invariance of these commutation relations yield new  
constraints on the elements of matrix $T$ which are consistent  
provided that  $k=\lambda_{12}=\mu_{12}=0$.
 \par
If instead we  choose
$(\xi^3)^2$ as third independent monomial, we obtain a different 
situation only if 
$\xi^1\xi^2=0$. In that case   $\xi^{\prime 1}\xi^{\prime 2}$ 
contains only one term of degree $0$, namely 
$T^1_{\, 3}T^2_{\, 3} (\xi^3)^2$ that cannot 
vanish. Hence this case is impossible.

We prove in the same way that  the dimension of $V_2$ cannot  be equal to 2.

{\bf Result}: 
If we assume that the differentials are 
related by six (or more) relations, the only possible set of relations 
is given by (\ref{13}) with $k=\lambda_{12}=\mu_{12}=0$. These 
relations will be denoted by $R^\prime _{\xi \xi}$.

The constraint of invariance of $R^\prime_{\xi\xi}$ has been explicitly 
studied in \cite{nous}. The consistency requirement     for  
the ordering of terms containing a product of three  $T^i_{\, j}$ 
has led to the condition  $ q=u^2$. As a consequence, the 
relations $R^\prime_{\xi\xi}$ become identical with
$R_{\xi\xi}$ given in (\ref{xi1}). 

%{\bf Result}
%$q=u^2$

\section{Quantum groups}
\subsection{A quantum group with seven generators}
A unique set of relations ($R_{xx},R_{\xi ,\xi}$) has been obtained. 
They are given by (\ref{3}, \ref{xi1}) with $q=u^2$   and 
define the $x$ and $\xi$ spaces respectively. The constraints 
of   invariance by homomorphism $\delta$ 
and consistency of the 
computations lead to  the following
 relations $R_{TT}$ between the elements 
$T^i_{\, j}$: 

\begin{equation}
\begin{array}{llll}
T^1_{\, 1} T^1_{\, 2} =u^{-2} \, T^1_{\, 2} T^1_{\, 1}  &    T^1_{\, 1} T^1_{\, 3} = u^{-1} \, T^1_{\, 3} T^1_{\, 1}   
&   T^1_{\, 1} T^2_{\, 2} =T^2_{\, 2} T^1_{\, 1}  & 
T^1_{\, 1} T^2_{\, 1}   =  u^2 \, T^2_{\, 1} T^1_{\, 1}      
\\
%[3mm]
  T^1_{\, 1} T^2_{\, 3} =u\, T^2_{\, 3} T^1_{\, 1}  &    T^1_{\, 1} T^3_{\, 3}  = T^3_{\, 3} T^1_{\, 1}    &
T^1_{\, 2} T^1_{\, 3} = u \, T^1_{\, 3} T^1_{\, 2}  &   T^1_{\, 2} T^2_{\, 2}   =  u^2 \, T^2_{\, 2} T^1_{\, 2}   
\\
%[3mm]  
T^1_{\, 2} T^2_{\, 1}  = u^4 \, T^2_{\, 1} T^1_{\, 2}   & 
T^1_{\, 2} T^2_{\, 3} = u^3 \, T^2_{\, 3} T^1_{\, 2}  &   
T^1_{\, 2} T^3_{\, 3} =u^2\,  T^3_{\, 3} T^1_{\, 2}  & 
T^1_{\, 3} T^2_{\, 2} = u\, T^2_{\, 2} T^1_{\, 3} 
 \\
% [3mm]
    T^1_{\, 3} T^2_{\, 1}  = u^3 \, T^2_{\, 1} T^1_{\, 3}  
&   
T^1_{\, 3} T^3_{\, 3}  = u\, T^3_{\, 3} T^1_{\, 3}  &  
T^2_{\, 2} T^2_{\, 1}  = u^2 \, T^2_{\, 1} T^2_{\, 2}  & T^2_{\, 2} T^2_{\, 3} = u \, T^2_{\, 3} T^2_{\, 2}   \\
%[3mm] 
    T^2_{\, 2} T^3_{\, 3}  = T^3_{\, 3} T^2_{\, 2}   & 
T^2_{\, 1} T^2_{\, 3}  = u^{-1} \, T^2_{\, 3} T^2_{\, 1}  
& T^2_{\, 1} T^3_{\, 3}  = u^{-2} \, T^3_{\, 3} T^2_{\, 1}   & 
T^2_{\, 3} T^3_{\, 3}  = u^{-1} \, T^3_{\, 3} T^2_{\, 3}  \\
%[3mm]
\multicolumn{4}{l}{
T^1_{\, 3} T^2_{\, 3} - u^2 \, T^2_{\, 3} T^1_{\, 3}  + s(T^1_{\, 1} T^2_{\, 2} - u^2 \, T^2_{\, 1} T^1_{\, 2}  -(T^3_{\, 3} )^2) = 0}   
 \end{array}
\label{22}
\end{equation}
The inverse matrix can now be completely determined and is given by :
\begin{equation}
T^{-1} = 
\left( 
\begin{array}{ccc}
T^2_{\, 2} T^3_{\, 3}  & -u^2 \, T^1_{\, 2} T^3_{\, 3}  & T^1_{\, 2} T^2_{\, 3}  -u\, T^1_{\, 3} T^2_{\, 2}  \\
-u^{-2} \, T^2_{\, 1} T^3_{\, 3}  &  T^1_{\, 1} T^3_{\, 3}  &  -u^{-2} \, T^1_{\, 1} T^2_{\, 3}  + u^{-3} \, T^1_{\, 3} T^2_{\, 1}  \\
0 & 0 &  T^1_{\, 1} T^2_{\, 2}  - u^{-2} \, T^1_{\, 2} T^2_{\, 1} 
\end{array}
\right)
\times D^{-1}
\label{19}
\end{equation}
\noindent where the determinant $D$ can be written as : 
\begin{equation}
D \equiv (T^1_{\, 1} T^2_{\, 2}  - u^{-2} \, T^1_{\, 2} T^2_{\, 1} )T^3_{\, 3}  
\label{20}
\end{equation}
Remark that the determinant is not central and its inverse $D^{-1}$ 
 must be added to the set of generators of the algebra. 
Its commutation relations $R_{T D^{-1}}$ are easily deduced from those of 
$D$ and read:
\begin{equation}
\begin{array}{lclcl}
D^{-1}T^1_{\, 1}  = T^1_{\, 1} D^{-1} & \hspace{3mm} & u^{-6}\, D^{-1}T^1_{\, 2}  = T^1_{\, 2} D^{-1} & \hspace{3mm} 
& u^{-3} \, D^{-1}T^1_{\, 3}  = T^1_{\, 3} D^{-1} \\
D^{-1}T^2_{\, 2}  = T^2_{\, 2} D^{-1} & \hspace{3mm} & D^{-1}T^2_{\, 1}  = u^{-6} \, T^2_{\, 1} D^{-1} & \hspace{3mm}  &
D^{-1}T^2_{\, 3}  = u^{-3} \, T^2_{\, 3} D^{-1} \\
D^{-1}T^3_{\, 3}  = T^3_{\, 3} D^{-1}
\end{array}
\label{23}
\end{equation}
With these definitions, it may be verified that 
$H_8 \equiv C<T,D^{-1}>/R_{TT}\cup R_{T D^{-1}}$ is a 
Hopf algebra with  co-product $\Delta$,  co-unit $\epsilon$ and antipode $S$  defined by :  
\begin{eqnarray}
\Delta(T) \equiv T \otimes T, \quad 
\Delta(D^{-1}) \equiv D^{-1} \otimes D^{-1}
\label{26} \\
\epsilon (T,D^{-1}) \equiv (I,1), \quad 
S(T)\equiv T^{-1}, \quad   S(D) \equiv D^{-1}
\label{28}
\end{eqnarray} 
Now we can apply  the usual method (\cite{wz},\cite{zumino}) 
to obtain the quadratic relations between the variables, the differentials 
and the derivatives \cite{nous}. They are given by:
\begin{eqnarray}
x^k \, \xi^l & 
\! \! \!  
= 
\! \! \!   
& \widehat{R}^{kl}_{mn} \, \xi^m  \, x^n \\
\partial_k \, \xi^l &
\! \! \!   
=
\! \! \!  
& \widehat{R}^{-1lm}_{\; \; \; \; kn}\, \xi^n  \, \partial_m \\
\partial_l \, x^k &
\! \! \!   
=
\! \! \!  
& \delta^k_l+\widehat{R}^{km}_{ln}\, x^n \, \partial_m
\end{eqnarray}

\subsection{A quantum group with nine generators}
 Once the  matrix ${\widehat{R}}$ has been explicitly computed, 
it is possible to introduce a new  quantum matrix $t$ with nine elements 
satisfying the relations $R_{tt}$ deduced from (\ref{66}): 
\begin{equation}
\widehat{R}^{ji}_{kl} \, 
t^k_{m} t^l_{n} = t^j_{l} t^i_{k} \, \widehat{R}^{lk}_{mn}
\label{666}
\end{equation}
The computation of the inverse $t^{-1}$ yields:
\begin{equation}
\left(\begin{array}{c c c}
t^2_2t^3_3-u \,t^2_3t^3_2&-u^2\, t^1_2t^3_3 + u^3\, 
t^1_3t^3_2&t^1_2t^2_3-u \,t^1_3t^2_2\\
-u^{-2}\, t^2_1t^3_3 +u^{-3} \,t^2_3t^3_1&t^1_1t^3_3-u^{-1} \,
t^1_3t^3_1&-u^{-2} \,t^1_1t^2_3+ u^{-3} \,t^1_3t^2_1\\
t^2_1t^3_2 -u^{-2} \,t^2_2t^3_1&-u^2 \,t^1_1t^3_2 +t^1_2t^3_1&
t^1_1t^2_2 -u^{-2}\, t^1_2t^2_1\\
\end{array}\right) d^{-1} 
\label{t1}
\end{equation}
with the determinant $d$ of $t$ equal to: 
$$d =t^1_1t^2_2t^3_3+t^1_3t^2_1t^3_2+u^{-3} \,t^1_2t^2_3t^3_1- u ^{-1}\,
t^1_1t^2_3t^3_2-u^{-2} \,t^1_2t^2_1t^3_3 -u^{-2}\, t^1_3t^2_2t^3_1.$$ 
It can be verified that $d$ is not a central element of $C<t>$ and 
therefore must be added  to this algebra. 
The commutation relations $R_{t d^{-1}} $ of $d^{-1}$ with the 
generators $t^i_j$ are:
\begin{equation}
\begin{array}{cccccc}
t^1_1 d^{-1}&=d^{-1} t^1_1,& t^1_2 d^{-1}&=u^{-6} \,d^{-1}t^1_2,& t^1_3 d^{-1}&= u^{-3} \,d^{-1}t^1_3,\\
t^2_2 d^{-1}&=d^{-1} t^2_2,& t^2_1 d^{-1}&= u^4 \,d^{-1}t^2_1,& t^2_3 d^{-1}&= u^{-3} \,d^{-1}t^2_3 ,\\
t^3_1 d^{-1} &= u^3 \,d^{-1} t^3_1,& t^3_2 d^{-1}&=u^{-3}\, d^{-1} t^3_2,& t^3_3 d^{-1}&= d^{-1}t^3_3.\\
\end{array}
\label{RD}
\end{equation}
In this manner, $H_{10} \equiv C<t,d^{-1}>/R_{tt}\cup R_{t d^{-1}}$ is endowed 
with a 
structure of Hopf algebra. \par 
Thus two Hopf algebras, $H_8$  and $H_{10}$, 
have been constructed. 
Both preserve the same differential calculus on the 
deformed oscillator algebra defined by $R_{xx}$ with  $q=u^2$.
In addition,  the construction ensures that $H_{10}$ contains 
   $H_8=C<T,D^{-1}>/R_{TT}\cup R_{T D^{-1}}$ as a Hopf subalgebra. 
\section{Conclusion}
We have been able to deform simultaneously the Weyl-Heisenberg 
algebra and its group of invariance (a subgroup of $GL(3)$). In addition, 
an invariant 
differential calculus     has been set up on the resulting quantum 
space. However, it must be stressed that the whole construction 
cannot be carried out for arbitrary values of the deformation 
parameters and that the final result depends only on one 
complex number $u$.
\par
The constraint on the values of the parameters can be removed when the 
requirement of invariance by a seven-generator quantum group is lifted. 
A purely algebraic approach \cite{zumino} can   be developed and the 
commutation relations  are then shown to be invariant by a 
quantum matrix belonging 
to a three-parameter deformation of $GL(3)$ \cite{irac}\cite{irac1}.
\par
The construction  performed in this paper 
has yielded two quantum groups and their associated 
 Hopf algebras  $H_{10}$, $H_8$, which have  
 ten  and eight generators respectively.
These algebras are original deformations of $GL(3)$ and of its subgroup $G$. 
Moreover, the smaller one $H_8$ is embedded in $H_{10}$ as a true 
Hopf subalgebra.
\vfill\eject


\begin{thebibliography}{99}
\bibitem{Woro} 
Woronowicz, S.L.(1989) Differential Calculus on Compact Matrix Pseudogroups 
(Quantum Groups) {\it Commun. Math. Phys.}, {\bf 122}, pp.~ 125-170.\par
\bibitem{Ber}
Bernard, D. (1990) Quantum Lie Algebras and Differential Calculus on 
Quantum Groups {\it Progress of Theoretical Physics Supplement }
{\bf 102} pp.~49-66.\par
\bibitem{wz}
Wess, J. and Zumino, B. (1990) Covariant Differential Calculus on the 
Quantum Hyperplane {\it Nucl. Phys. (Proc. Suppl.)}, {\bf B 18}, pp.~301-313. \par
\bibitem{zumino}
Zumino, G. (1992) Differential calculus on quantum spaces and quantum 
groups, in {\it Group Theoretical Methods in Physics}, M.A. del Olmo, 
M.Santander and J.Mateos Guilarte (Eds.), Anales de F\'{\i}sica 
Monografias, CIEMAT (Spain), pp.~41-59.\par
\bibitem{rtf}

Reshetikhin, N.Yu., Takhtadzhyan, L.A. and Faddeev, L.D. (1990) 
Quantization of Lie groups and Lie algebras, {\it Leningrad Math. J.}, 
{\bf 1}, pp.~193-225.\par

 \bibitem{nous}
Bertrand, J. and Irac-Astaud, M. (1995) Invariant differential calculus on 
a deformation of the Weyl-Heisenberg algebra, {\it Modern Group Theoretical 
Methods in Physics}, p.~37, Kluwer Academic Publishers.
\bibitem{irac}
Irac-Astaud, M.(1996) Differential calculus on a three-parameter 
oscillator algebra,  to be published in {\it Reviews in Mathematical Physics}
\bibitem{irac1}
Irac-Astaud, M.(1996) A three-parameter deformation of the Weyl-Heisenberg 
algebra : differential calculus and invariance, to be published in the 
Proceedings of the 5th international Colloquium 
"Quantum Groups and Integrable systems" of  Prague ,  
{\it Czechoslovak Journal of Physics}.

\end{thebibliography}
\end{document}